\documentclass[12pt]{article}

\usepackage{amsmath}
\usepackage{amsfonts}
\usepackage{amssymb}
\usepackage{theorem}
\usepackage{graphicx}
\usepackage[dcucite]{harvard}

\theorembodyfont{\it}
\newtheorem{theorem}{Theorem}[section]

\newtheorem{lemma}[theorem]{Lemma}

\newtheorem{star-theorem}{Theorem}

\theorembodyfont{\rm}
\newtheorem{definition}[theorem]{Definition}
\newtheorem{notation}[theorem]{Notation}

\def\proof{\noindent{\bf Proof: }}
\def\endproof{\hfill$\square$}

\renewenvironment{enumerate}
{
\begin{list}{\makebox[1.0cm]{\emph{(\roman{enumi})}}}
{
\usecounter{enumi}
\setlength{\topsep}{0pt} \setlength{\partopsep}{0pt} \setlength{\parsep}{0pt}
\setlength{\itemsep}{0pt} \setlength{\labelsep}{0in}
\setlength{\labelwidth}{1.0cm} \setlength{\leftmargin}{1.0cm} } } {
\end{list}
}

\renewcommand{\harvardand}{and}

\begin{document}

\thispagestyle{empty}

\title{Nilpotentization of the Kinematics of the N-Trailer
       System at Singular Points and Motion Planning
       Through the Singular Locus}

\author{William {\sc Pasillas-L\'epine}\thanks{Now at Laboratoire des signaux et syst\`emes,
CNRS -- Sup\'elec. Plateau de Moulon, 3 rue Joliot-Curie. 91 192 Gif-sur-Yvette
Cedex, France. E-mail: \texttt{pasillas@lss.supelec.fr}} \hbox{ and }Witold
{\sc Respondek} \bigskip \\ Institut national des sciences appliqu\'ees de
Rouen\\ D\'epartement g\'enie math\'ematique\\ Place \'Emile Blondel \\ 76 131
Mont Saint Aignan Cedex, France \bigskip \\ Tel: (+33) 02 35 52 84 32 ; Fax:
(+33) 02 35 52 83 32 \\ E-Mail: {\tt wresp@lmi.insa-rouen.fr}}

\date{March, 2000}

\enlargethispage{2cm}

\maketitle

\begin{abstract}
We propose in this paper a constructive procedure that transforms locally, even
at singular configurations, the kinematics of a car towing trailers into
Kumpera-Ruiz normal form. This construction converts the nonholonomic motion
planning problem into an algebraic problem (the resolution of a system of
polynomial equations), which we illustrate by steering the two-trailer system
in a neighborhood of singular configurations. We show also that the $n$-trailer
system is a universal local model for all Goursat structures and that all
Goursat structures are locally nilpotentizable.

\medskip

\noindent {\bf Keywords:} Nonholonomic control systems, feedback equivalence,
Kumpera-Ruiz normal forms, Goursat normal form, car with trailers, nilpotent
Lie algebras, nonholonomic motion planning.

\end{abstract}

\newpage

\section*{Introduction}

The nonholonomic motion planning problem for a car towing trailers has been one
of the most widely studied problems in nonlinear control theory during the last
ten years. The controllability of this system was proved
by~\citeasnoun{laumond-trailer}.
%both for regular and singular configurations.
Since then many important results have been obtained, including solutions of
the motion planning and stabilization problems --- see the papers
\cite{fliess-intro-flat,fliess-trailer,jiang-nijmeijer,samson-chained,sordalen-trailer,teel-murray-walsh,tilbury-murray-sastry},
the book \cite{laumond-book}, and references given there --- with the main
emphasis put on regular configurations. For singular configurations much less
results are available
\citeaffixed{cheaito-mormul,jean-trailer,pasillas-respondek-nolcos,vendittelli-laumond-oriolo}{see,
however,}, although it is clear that for some steering situations, like turning
and going back in a narrow street, it is necessary to cross the singular locus.
The aim of this paper is to show that theoretical results of
\cite{cheaito-mormul,kumpera-ruiz,pasillas-respondek-goursat} can be applied to
this problem. We illustrate the theory with a concrete example: the two-trailer
system.

The paper is organized as follows. In
Section~\ref{sec-ieee-goursat-and-kumpera-ruiz} we introduce Goursat structures
and define Kumpera-Ruiz normal forms for Goursat structures. In
Section~\ref{sec-nilpotentization} we show that any Kumpera-Ruiz normal form is
defined by a pair of vector fields that generate a nilpotent Lie algebra thus
proving that any Goursat structure is locally feedback nilpotentizable. In
Section~\ref{sec-ieee-n-trailer} we recall the nonholonomic model of a car
towing trailers. Then we show how to transform the $n$-trailer system into a
Kumpera-Ruiz normal form and we prove a surprising result: any Goursat
structure is locally equivalent to the $n$-trailer system around a well chosen
point of its configuration space. An alternative proof of this fact has been
proposed by \citeasnoun{montgomery-zhitomirskii}. In
Section~\ref{sec-motion-planning}, we briefly discuss how transforming the
$n$-trailer system to a Kumpera-Ruiz normal form simplifies the motion planning
problem by reducing it to an algebraic problem of solving a system of
polynomial equations. We illustrate this approach in
Section~\ref{sec-two-trailer} by applying it to the motion planning problem
through singular configurations for the two-trailer system. Preliminary results
have been published in \cite{pasillas-respondek-cdc}.

\noindent {\bf Acknowledgments:} The authors would like to thank Henk Nijmeijer
for interesting comments on the paper.

\section{Goursat Structures and Kumpera-Ruiz Normal Forms}

\label{sec-ieee-goursat-and-kumpera-ruiz}

A rank-$k$ \emph{distribution} $\mathcal D$ on a smooth manifold $M$ is a map
that assigns smoothly to each point $p$ in $M$ a linear subspace
$\mathcal{D}(p) \subset T_pM$ of dimension $k$. Such a field of $k$-planes  is
locally spanned by $k$ pointwise linearly independent smooth vector fields
$f_{1},\ldots,f_{k}$, which will be denoted by
$\mathcal{D}=(f_{1},\ldots,f_{k})$.

Two distributions $\mathcal{D}$ and $\tilde{\mathcal{D}}$ defined on two
manifolds~$M$ and $\tilde M$, respectively, are \emph{equivalent} if there
exists a smooth diffeomorphism $\varphi$ between $M$ and $\tilde M$ such that
$(\varphi_{*}\mathcal{D})(\tilde p)=\tilde{\mathcal{D}}(\tilde p)$, for each
point $\tilde p$ in $\tilde M$. Clearly, two distributions $\mathcal{D}$ and
$\tilde{\mathcal{D}}$ are locally equivalent if and only if the corresponding
control systems $$ \dot{x}={\textstyle \sum\limits_{i=1}^{k}}f_{i}(x)\,u_{i}
\quad \text{ and } \quad \dot{\tilde{x}}={\textstyle
\sum\limits_{i=1}^{k}}\tilde{f}_{i}(\tilde{x})\,\tilde{u}_{i} $$ are locally
equivalent via a state static feedback.

The \emph{derived flag} of a distribution $\mathcal{D}$ is the sequence of
modules $\mathcal{D}^{(0)}\subset\mathcal{D}^{(1)}\subset\cdots$ defined
inductively by $\mathcal{D}^{(0)}=\mathcal{D}$ and
$\mathcal{D}^{(i+1)}=\mathcal{D}^{(i)}+[\mathcal{D}^{(i)},\mathcal{D}^{(i)}]$,
for $i\geq0$.

\begin{definition}
A \emph{Goursat structure} on a manifold~$M$ of dimension $n\geq3$ is a
rank-two distribution~$\mathcal{D}$ such that, for $0\leq i\leq n-2$, the
elements of its derived flag satisfy $\dim\mathcal{D}^{(i)}(p)=i+2$, for any
point$~p$ in~$M$.
\end{definition}

Any Goursat structure on a manifold $M$ of dimension~$n$ is equivalent, in a
small enough neighborhood of any point of an open and dense subset of $M$, to
the distribution spanned by
\[
\left(
\begin{array}[c]{c}
\tfrac\partial{\partial x_{n}}
\end{array}
,
\begin{array}
[c]{c}%
x_{n}\tfrac\partial{\partial x_{n-1}}+\cdots+x_{3}\tfrac\partial{\partial
x_{2}}+\tfrac\partial{\partial x_{1}}%
\end{array}
\right),
\]
which is called \emph{Goursat normal form} or \emph{chained form}. We refer the
reader to
\cite{bryant-chern-gardner-goldschmidt-griffiths,kumpera-ruiz,martin-rouchon-driftless,murray-nilpotent}
for additional information about this normal form. If at a given point a
Goursat structure can be converted into Goursat normal form then this point is
called \emph{regular}; otherwise, it is called \emph{singular}. An elegant
characterization of regular points, based on the concept of growth vector, can
be found in the work of~\citeasnoun{murray-nilpotent}.

\medskip

For Goursat structures, the first who observed the existence of singular points
were~\citeasnoun{giaro-kumpera-ruiz}. This initial observation has
led~\citeasnoun{kumpera-ruiz} to write their innovative article, where they
introduced the normal forms that we consider in this section.

We start with the Pfaff-Darboux and Engel normal forms, given respectively
on~$\mathbb{R}^{3}$ and~$\mathbb{R}^{4}$ by the pairs of vector fields
$\kappa^{3}=(\kappa_{1}^{3},\kappa_{2}^{3})$ and $\kappa^{4}=(\kappa_{1}^{4},
\kappa_{2}^{4})$, defined respectively by
\[
\begin{array}[c]{l}
\kappa_{1}^{3}=\tfrac\partial{\partial x_{3}}\\
\kappa_{2}^{3}=x_{3}\tfrac\partial{\partial x_{2}}+\tfrac\partial{\partial
x_{1}}%
\end{array}
\qquad\qquad
\]
and
\[
\begin{array}[c]{l}
\kappa_{1}^{4}=\tfrac\partial{\partial x_{4}}\\
\kappa_{2}^{4}=x_{4}\tfrac\partial{\partial x_{3}}+x_{3}\tfrac\partial
{\partial x_{2}}+\tfrac\partial{\partial x_{1}}.
\end{array}
\]
Loosely speaking, we can write
\[
\begin{array}[c]{l}
\kappa_{1}^{4}=\tfrac\partial{\partial x_{4}}\\
\kappa_{2}^{4}=x_{4}\kappa_{1}^{3}+\kappa_{2}^{3}.
\end{array}
\qquad\qquad\qquad
\]
In order to make this precise we will adopt the following natural convention.
Consider a vector field $$f^{n-1}=%
%TCIMACRO{\tsum _{i=1}^{n-1}}%
%BeginExpansion
{\textstyle\sum\limits_{i=1}^{n-1}}
%EndExpansion
f_{i}^{n-1}(x_{1},\ldots,x_{n-1})\tfrac\partial{\partial x_{i}}$$ on
$\mathbb{R}^{n-1}$ equipped with coordinates $(x_{1},\ldots,x_{n-1})$. We can
lift $f^{n-1}$ to a vector field, denoted also by $f^{n-1} $, on
$\mathbb{R}^{n}$ equipped with coordinates $(x_{1},\ldots,x_{n-1},x_{n})$ by
taking
$$f^{n-1}=%
%TCIMACRO{\tsum _{i=1}^{n-1}}%
%BeginExpansion
{\textstyle\sum\limits_{i=1}^{n-1}}
%EndExpansion
f_{i}^{n-1}(x_{1},\ldots,x_{n-1})\tfrac\partial{\partial x_{i}}+0\cdot\tfrac
\partial{\partial x_{n}}.$$ That is, we lift $f^{n-1}$ by translating it along
the $x_{n}$-direction.

\begin{notation}
\label{not-ieee-lift} From now on, in any expression of the form $\kappa
_{2}^{n}=\alpha(x)\kappa_{1}^{n-1}+\beta(x)\kappa_{2}^{n-1}$, the vector fields
$\kappa_{1}^{n-1}$ and $\kappa_{2}^{n-1} $ should be considered as the above
defined lifts of $\kappa_{1}^{n-1}$ and $\kappa_{2}^{n-1}$, respectively.
\end{notation}

Let $\kappa^{n-1}=(\kappa_{1}^{n-1},\kappa_{2}^{n-1})$ denote a pair of vector
fields on $\mathbb{R}^{n-1}$. A \emph{regular prolongation,\ with parameter}
$c_{n}$, of~$\kappa^{n-1}$, denoted by $\kappa^{n}=R_{c_{n}}(\kappa^{n-1})$, is
a pair of vector fields $\kappa^{n}=(\kappa_{1}^{n},\kappa_{2}^{n})$ defined on
$\mathbb{R}^{n}$ by
\begin{equation}%
\begin{array}
[c]{l}%
\kappa_{1}^{n}=\tfrac\partial{\partial x_{n}}\\
\kappa_{2}^{n}=(x_{n}+c_{n})\kappa_{1}^{n-1}+\kappa_{2}^{n-1},
\end{array}
\label{ieee-regular-prolongation}%
\end{equation}
where $c_{n}$ belongs to $\mathbb{R}$. The \emph{singular prolongation}
of~$\kappa^{n-1}$, denoted by $\kappa^{n}=S(\kappa^{n-1})$, is the pair of
vector
fields $\kappa^{n}=(\kappa_{1}^{n},\kappa_{2}^{n})$ defined on $\mathbb{R}%
^{n}$ by
\begin{equation}%
\begin{array}
[c]{l}%
\kappa_{1}^{n}=\frac\partial{\partial x_{n}}\\
\kappa_{2}^{n}=\kappa_{1}^{n-1}+x_{n}\kappa_{2}^{n-1}.
\end{array}
\qquad\quad\label{ieee-singular-prolongation}%
\end{equation}

\begin{definition}
For $n\geq3$, a pair of vector fields $\kappa^{n}$ on $\mathbb{R}^{n}$ is
called a \emph{Kumpera-Ruiz normal form} if it is given by
%it can be obtained by a sequence
%of prolongations from the Pfaff-Darboux normal form. In other words, we have
$\kappa^{n}=\sigma_{n-3}\circ\cdots\circ\sigma_1(\kappa^{3})$, where each
$\sigma_i$, for $1\leq i\leq n-3$, equals either $S$ or $R_{c_i}$, for some
real constants~$c_i$.
\end{definition}

Note that our definition differs slightly from that of
\citeasnoun{cheaito-mormul-pasillas-respondek}. Firstly, we do not ask the
coordinates
%$x:\mathbb{R}^{n}\rightarrow\mathbb{R}^{n}$
to satisfy $x(p)=0$, where $p$ is the point around which we work. Secondly, we
consider the models
$S(\kappa^{3})$ and $R_{c}(\kappa^{3})$, which are equivalent to $R_{0}%
(\kappa^{3})$, as being Kumpera-Ruiz normal forms. The following result of
\citeasnoun{kumpera-ruiz}
\citeaffixed{cheaito-mormul,cheaito-mormul-pasillas-respondek,montgomery-zhitomirskii,pasillas-respondek-goursat}{see
also} shows clearly the interest of their normal forms.

\begin{theorem}[Kumpera-Ruiz]
\label{ieee-thm-kumpera-ruiz}Any Goursat structure is locally equivalent to a
distribution spanned by a Kumpera-Ruiz normal form.
\end{theorem}

This result implies that locally, even at singular points, Goursat structures
do not have functional invariants. This make them precious but rare and
distinguish them from generic rank-two distributions on $n$-manifolds, which do
have local functional invariants when $n\geq5$. Moreover, as we will see in the
next Section, this result implies that any Goursat structure is locally
feedback nilpotentizable.

\section{Nilpotentization}

\label{sec-nilpotentization}

Let us recall the following standard concepts
\citeaffixed{fulton-representation-theory}{see e.g.}. Let~$\mathfrak{g}$ be a
Lie algebra. A \emph{Lie subalgebra} of $\mathfrak{g}$ is a
linear subspace $\mathfrak{h}\subset\mathfrak{g}$ such that $[\mathfrak{h},\mathfrak{h}%
]\subset\mathfrak{h}$. An \emph{ideal} of $\mathfrak{g}$ is a Lie subalgebra
$\mathfrak{i}$ such that $[\mathfrak{i},\mathfrak{g}]\subset\mathfrak{i}$. The
\emph{lower central series} of a Lie algebra $\mathfrak{g}$ is the sequence of
ideals
$\mathfrak{D}_{0}(\mathfrak{g})\supset\mathfrak{D}_{1}(\mathfrak{g})\supset\cdots$
defined by
\[
\mathfrak{D}_{0}(\mathfrak{g})=\mathfrak{g}\text{\quad and\quad}\mathfrak{D}_{k}%
(\mathfrak{g})=[\mathfrak{g},\mathfrak{D}_{k-1}(\mathfrak{g})]\text{,\quad for
}k\geq1\text{. }
\]
In other words, for a fixed $k\geq0$ the ideal $\mathfrak{D}_{k}(\mathfrak{g})$
is the subspace generated by the elements of $\mathfrak{g}$ that can be
expressed as a left-iterated Lie bracket of the form
$[e_{0},[e_{1},\ldots[e_{l-1},e_{l}]\ldots ]]$, where $l\geq k$ and
$e_{0},\ldots,e_{l}$ are elements of $\mathfrak{g}$. Observe that unlike the
Lie and derived flags, which increase, the lower central series decreases.

A Lie algebra $\mathfrak{g}$ is \emph{nilpotent} if there exists some integer
$r$ such that $\mathfrak{D}_{r}(\mathfrak{g})=0$; the smallest such integer is
called the \emph{nilindex} of $\mathfrak{g}$. It is clear that a Lie algebra
$\mathfrak{g}$ is nilpotent if and only if there exists some integer $r$ such
that, for every sequence $e_{0},\ldots,e_{r}$ of elements of $\mathfrak{g}$, we
have
\[
\lbrack e_{0},[e_{1},\ldots[e_{r-1},e_{r}]\ldots]]=0.
\]
Denote $\mathrm{ad}_{\mathrm{\,}e_{0}}(e_{1})=[e_{0},e_{1}]$, for all $e_{0}$
and $e_{1}$ in $\mathfrak{g}$ . For each element $e_{0}$ in $\mathfrak{g}$, the
map $\mathrm{ad\,}_{e_{0}}$ is linear. The following result is standard
\citeaffixed{fulton-representation-theory}{see e.g.}.

\begin{lemma}[Engel]
\label{lem-engel}A finite dimensional Lie algebra $\mathfrak{g}$ is nilpotent
if and only if $\mathrm{ad\,}_{e_0}$ is nilpotent for each element $e_0$ in
$\mathfrak{g}$.
\end{lemma}

For any subspace $V$ of a given Lie algebra $\mathfrak{g}$, not necessarily of
finite dimension, define the sequence of subspaces $V^{(0)}\subset
V^{(1)}\subset\cdots$ by
\[
V^{(0)}=V\text{\quad and\quad}V^{(k+1)}=V^{(k)}+[V^{(k)},V^{(k)}].
\]
The Lie algebra \emph{generated} by a subspace $V\subset\mathfrak{g}$ is the
infinite sum
\[
\mathcal{L}(V)=V^{(0)}+V^{(1)}+\cdots+V^{(k)}+\cdots,
\]
which is clearly a Lie subalgebra of $\mathfrak{g}$.

A distribution is said to be $\emph{nilpotentizable}$ if we can chose a family
of vector fields that span the distribution and generate a nilpotent Lie
algebra of finite dimension. The class of nilpotentizable distributions is
particularly important in control theory, because for them, for instance, a
general motion planning algorithm exists
\citeaffixed{hermes-nilpotent,laferriere-sussmann,murray-nilpotent}{see e.g.}.
The following result states that any nonholonomic control system with two
controls that generates a Goursat structure is locally feedback equivalent to a
system whose Lie algebra is finite dimensional and nilpotent.

\begin{theorem}
\label{thm-nilpotentizable}Goursat structures are locally nilpotentizable. In
fact, any Kumpera-Ruiz normal form generates a nilpotent Lie algebra that has
finite dimension.
\end{theorem}

\proof By Theorem~\ref{ieee-thm-kumpera-ruiz}, any Goursat structure on a
manifold of dimension $n$
is locally equivalent to a Kumpera-Ruiz normal form $(\kappa_{1}^{n}%
,\kappa_{2}^{n})$ centered at zero. It thus suffices to prove that
$(\kappa_{1}^{n},\kappa_{2}^{n}) $ generates a nilpotent Lie algebra.

Let $(\kappa_{1}^{n},\kappa_{2}^{n})$ be any Kumpera-Ruiz normal form on
$\mathbb{R}^{n}$. Denote by $\mathfrak{h}_{n}$ the Lie algebra generated by
$\kappa_{1}^{n}$ and $\kappa_{2}^{n}$. The main argument of the proof is to
show that~$\mathfrak{h}_{n}$ is a Lie subalgebra of a nilpotent Lie
algebra~$\mathfrak{g}_{n}$ of finite dimension. Thus $\mathfrak{h}_{n}$ is
nilpotent, since any subalgebra of a nilpotent Lie algebra is itself nilpotent.

Let us proceed by induction on $n$. For $n=3$ it is clear that $\mathfrak
{h}_{3}$ is nilpotent, since~$\mathfrak{h}_{3}$ is the three-dimensional
Heisenberg algebra. Assume that $\mathfrak{h}_{n-1}$ is nilpotent and of finite
dimension. We will prove that this assumption implies that $\mathfrak{h}_{n}$
is nilpotent and of finite dimension.

Put $V=\operatorname{span}_{\mathbb{R}}(\kappa_{1}^{n},\kappa_{2}^{n})$. From
the definition of $V^{(0)}$ and $V^{(1)}$, and the definition of Kumpera-Ruiz
normal forms, we have
\[
V^{(0)}=\operatorname{span}_{\mathbb{R}}(\kappa_{1}^{n},\kappa_{2}%
^{n})\text{\quad and\quad}V^{(1)}\subset\operatorname{span}_{\mathbb{R}%
}(\kappa_{1}^{n},\kappa_{2}^{n},\mathfrak{h}_{n-1}).
\]
Define the sequence $E^{(i)}$, for $i\geq0$, of subspaces of the Lie algebra
of all polynomial vector fields on~$\mathbb{R}^n$, by taking $E^{(0)}=\operatorname{span}_{\mathbb{R}%
}(\kappa_{1}^{n},\kappa_{2}^{n})$ and
\begin{equation}
E^{(i)}=E^{(i-1)}+P_{2^{i-1}}[x_{n}]\mathfrak{D}_{i-1}(\mathfrak{h}%
_{n-1}),\text{\quad for }i\geq1,\label{nilpotent-sequence}%
\end{equation}
where $P_{k}[x_{n}]$ denotes the vector space of polynomials with real
coefficients of the variable $x_{n}$ that have degree $d\leq k$, and
$P_{2^{i-1}}[x_{n}]\mathfrak{D}_{i-1}(\mathfrak{h}_{n-1})$ denotes the subspace
obtained by taking linear combinations, with coefficients in
$P_{2^{i-1}}[x_{n}]$, of vectors fields in $\mathfrak{D}_{i-1}(\mathfrak
{h}_{n-1})$.

Denote by $r$ the nilindex of $\mathfrak{h}_{n-1}$ and put
$\mathfrak{g}_{n}=E^{(r)}$. The rest of the proof is a direct consequence of
the following Lemma, which shows that $\mathfrak{h}_{n}$, the Lie algebra
generated by $V^{(0)}$, is contained in the nilpotent Lie
algebra~$\mathfrak{g}_{n}$, that has finite dimension
\endproof

\begin{lemma}
We have the following properties of $V^{(i)}$ and $E^{(i)}$:

\begin{enumerate}
\item $E^{(i)}+[E^{(i)},E^{(i)}]\subset E^{(i+1)};$

\item $V^{(i)}\subset E^{(i)};$

\item $E^{(r+1)}=E^{(r)};$

\item $E^{(r)}$ is a Lie algebra of finite dimension;

\item $E^{(r)}$ is nilpotent.
\end{enumerate}
\end{lemma}%

\proof Let us denote $W_{0}^{n}=\operatorname{span}_{\mathbb{R}}(\kappa_{1}%
^{n},\kappa_{2}^{n})$ and $$W_{i}^{n}=P_{2^{i-1}}[x_{n}]\mathfrak{D}_{i-1}%
(\mathfrak{h}_{n-1})\mbox{, for }i\geq1.$$ \emph{First Item}. This Item can be
proved by induction on $i$. We clearly have $E^{(0)}+[E^{(0)},E^{(0)}]\subset
E^{(1)}$. Assume that Item (i) is true up to $i$. We are going to prove that
then it is also true for $i+1$. Indeed, we have
\[
E^{(i)}=W_{0}^{n}+W_{1}^{n}+\cdots+W_{i}^{n}.
\]
Since, by the induction assumption, $E^{(i)}+[E^{(i)},E^{(i)}]\subset
E^{(i+1)}$, what remains to prove is that $[W_{j}^{n},W_{i+1}^{n}]\subset
E^{(i+2)}$, for $0 \leq j \leq i+1$. But this last relation results directly
from the definitions of $\mathfrak{D}_{i+1}(\mathfrak{h}_{n-1})$ and
$P_{2^{i+1}}[x_{n}]$, and from standard properties of the Lie bracket.

\emph{Second Item}. This Item is a direct consequence of Item (i) and the
relation $V^{(0)}=E^{(0)}$.

\emph{Third Item}. Since the nilindex of $\mathfrak{h}_{n-1}$ is $r$, the
sequence (\ref{nilpotent-sequence}) stabilizes at
\[
E^{(r)}=E^{(r-1)}+P_{2^{r-1}}[x_{n}]\mathfrak{D}_{r-1}(\mathfrak{h}%
_{n-1})\text{.}
\]
Indeed, we have $E^{(r+1)}=E^{(r)}+P_{2^{r}}[x_{n}]\mathfrak{D}_{r}(\mathfrak
{h}_{n-1})$ and thus $E^{(r+1)}=E^{(r)}$, since $\mathfrak{D}_{r}(\mathfrak
{h}_{n-1})=0$.

\emph{Fourth Item}. Items (i) and (iii) imply $E^{(r)}+[E^{(r)},E^{(r)}]\subset
E^{(r)}$, that is $E^{(r)}$ is closed under Lie brackets. But since, by
construction, to the finite dimensional vector space $E^{(i)}$ we add (at each
step) a finite dimensional vector space (recall that $\mathfrak{h}_{n-1}$ has
finite dimension), the dimension of $E^{(r)}$ is finite.

\emph{Fift Item}. It is a direct consequence of Engel's Lemma that $E^{(r)}$ is
nilpotent, because $\mathrm{ad\,}_{e_{0}}$ is nilpotent for any $e_{0}$ in
$E^{(r)}$. Indeed, observe that  for any $w_i \in W_i^n$ and $w_j \in W_j^n$,
such that $i \leq j$, we have $[w_i,w_j] \in W_{j+1}^n$.

\endproof

\section{The N-Trailer System}

\label{sec-ieee-n-trailer}

The kinematical model for a unicycle-like mobile robot towing $n$ trailers such
that the tow hook of each trailer is located at the center of its unique axle
is usually called, in control theory, the $n$-trailer system --- see
\cite{fliess-intro-flat},  \cite{jean-trailer}, \cite{laumond-trailer},
\cite{laumond-book}, \cite{fliess-trailer}, \cite{samson-chained},
\cite{sordalen-trailer}, \cite{teel-murray-walsh},
\cite{tilbury-murray-sastry}, and references therein. For simplicity, we will
assume that the distances between any two consecutive trailers are equal.

We give here an inductive definition of the $n$-trailer. This definition
already appears in~\cite{jean-trailer} and reminds the one given in the
previous section for Kumpera-Ruiz normal forms. To start with, consider the
pair of vector fields $(\tau_{1}^{0},\tau_{2}^{0})$ on $\mathbb{R}^{2}\times
S^{1}$, equipped with coordinates $(\xi_{1},\xi_{2},\theta_{0})$, that describe
the kinematics of the unicycle-like mobile robot towing no trailers:
\[
\begin{array}[c]{l}
\tau_{1}^{0}=\tfrac\partial{\partial\theta_{0}}\\
\tau_{2}^{0}=\sin(\theta_{0})\tfrac\partial{\partial\xi_{2}}+\cos(\theta
_{0})\tfrac\partial{\partial\xi_{1}}.
\end{array}
\]
The \emph{$n$-trailer system} is defined by applying successively a sequence of
prolongations to this mobile robot. In order to do this, suppose that a pair of
vector fields $\tau^{n-1}=(\tau_{1}^{n-1},\tau_{2}^{n-1})$ on
$\mathbb{R}^{2}\times(S^{1})^{n}$ associated to the mobile robot towing $n-1$
trailers has been defined. The pair of vector fields $\tau^{n}=(\tau_{1}%
^{n},\tau_{2}^{n})$ on $\mathbb{R}^{2}\times(S^{1})^{n+1}$ defining the
$n$-trailer system is given by
\[
\begin{array}[c]{l}
\tau_{1}^{n}=\tfrac\partial{\partial\theta_{n}}\\
\tau_{2}^{n}=\sin(\theta_{n}-\theta_{n-1})\tau_{1}^{n-1}+\cos(\theta
_{n}-\theta_{n-1})\tau_{2}^{n-1},
\end{array}
\]
where the coordinates $\xi_{1}$ and $\xi_{2}$ represent the position of the
last trailer, while the coordinates $\theta_{0},\ldots,\theta_{n}$ represent,
starting from the last trailer, the angles between each trailer's axle and the
$\xi_{1}$-axis. Observe that this definition should be understood in the sense
of Notation~\ref{not-ieee-lift}. Mechanically, to prolongate the $n$-trailer
means to add one more trailer to the system. It is easy to check that the
$n$-trailer system is a Goursat structure.

\medskip

A direct consequence of Kumpera-Ruiz's theorem is that, in a small enough
neighborhood of any point of its configuration space, in particular at any
singular configuration, the $n$-trailer can be converted into Kumpera-Ruiz
normal form. One of the main results of this paper is to describe this
conversion explicitly. For regular configurations, our result gives the
transformations proposed in \cite{sordalen-trailer} and
\cite{tilbury-murray-sastry}; for singular configurations, a new kind of
transformations is obtained.

The $n$-trailer system can also be written as the control system
\begin{equation}
\dot\zeta=\tau_{1}^{n}(\zeta)\,v_{1}+\tau_{2}^{n}(\zeta)\,v_{2},
\label{trailer-control-system}%
\end{equation}
where $\zeta=(\xi_{1},\xi_{2},\theta_{0},...,\theta_{n})$.
%denotes local coordninates on $\mathbb{R}^{2}\times(S^{1})^{n+1}$.
Recall that a feedback transformation of (\ref{trailer-control-system}) is a
change of controls of the form
\begin{equation}%
\begin{array}
[c]{l}%
u_{1}=\nu(\zeta)\,v_{1}+\eta(\zeta)\,v_{2}\\
u_{2}=\lambda(\zeta)\,v_{1}+\mu(\zeta)\,v_{2},
\end{array}
\label{feedback-transformation}%
\end{equation}
where the smooth functions $\nu$, $\lambda$, $\eta$, and $\mu$ are such that
$(\nu\mu-\lambda\eta)(\cdot)\neq0$.

Fix a point $p$ of $\mathbb{R}^{2}\times(S^{1})^{n+1}$ given in $\zeta
$-coordinates by
\[
\zeta(p)=\zeta^p=(\xi_{1}^{p},\xi_{2}^{p},\theta_{0}^{p},...,\theta_{n}^{p}).
\]
In order to convert, locally at~$p$, the $n$-trailer into a Kumpera-Ruiz normal
form we look for a local change of coordinates
$$(x_{1},\ldots,x_{n+3})=\phi^{n}(\xi_{1},\xi_{2},\theta_{0},\ldots,\theta
_{n})$$ and a local triangular feedback transformation of the form
$u_{1}=\nu_{n}(\zeta)\,v_{1}+\eta_{n}(\zeta)\,v_{2}$ and
$u_{2}=\mu_{n}(\zeta)\,v_{2}$ that bring~(\ref{trailer-control-system}) into
\begin{equation}
\dot x=\kappa_{1}^{n}(x)\,u_{1}+\kappa_{2}^{n}(x)\,u_{2},
\label{kumpera-ruiz-control-system}
\end{equation}
where $x=(x_{1},\ldots,x_{n+3})$. In other words, we ask the change of
coordinates and the feedback transformation to satisfy
\begin{equation}
\begin{array}
[c]{l}%
\kappa_{1}^{n+3}=\phi_{*}^{n}(\widehat{\nu}_{n}\tau_{1}^{n})\\
\kappa_{2}^{n+3}=\phi_{*}^{n}(\widehat{\eta}_{n}\tau_{1}^{n}+\widehat{\mu}%
_{n}\tau_{2}^{n})
\end{array}
\quad\quad\qquad\label{kumpera-as-a-trailer}%
\end{equation}
or, equivalently,
\begin{equation}%
\begin{array}[c]{lll}
\phi_{*}^{n}(\tau_{1}^{n}) & = & (\nu_{n}\circ\psi^{n})\,\kappa_{1}^{n+3}\\
\phi_{*}^{n}(\tau_{2}^{n}) & = & (\eta_{n}\circ\psi^{n})\,\kappa_{1}^{n+3}
+(\mu_{n}\circ\psi^{n})\,\kappa_{2}^{n+3},
\end{array}
\label{ieee-trailer-as-kumpera}%
\end{equation}
where $\psi^{n}=(\phi^{n})^{-1}$ denotes the inverse of the local
diffeomorphism $\phi^{n}$, both
%$\nu_{n}(p)\neq0$ and $\mu_{n}(p)\neq0$
$\nu_{n}(\zeta^p)\neq0$ and $\mu_{n}(\zeta^p)\neq0$, and the inverse feedback
transformation
$v_{1}=\widehat{\nu}_{n}(\zeta)\,u_{1}+\widehat{\eta}_{n}(\zeta)\,u_{2}$ and
$v_{2}=\widehat{\mu}_{n}(\zeta)\,u_{2}$ is obviously given by
\begin{equation}
\label{inverse-feedback-transformation} \widehat{\nu}_{n}=\frac1{\nu_{n}},
\text{\quad}\widehat{\eta}_{n}=-\frac{\eta_{n}}{\mu_{n}\nu_{n}}, \text{\quad
and\quad}\widehat{\mu}_{n}=\frac1{\mu_{n}}.
\end{equation}
Observe that we not demand the $x$-coordinates to be centered at $p$, and thus
the point $x(p)=(\phi^{n}\circ\zeta)(p)$ will be, in general, different from
zero.

To start our construction, take $x_{1}=\xi_{2}$ and $x_{2}=\xi_{1}$. If
$\theta_{0}^{p}\neq\pm\pi/2\mod2\pi$ take $x_{3}=\tan(\theta_{0})$, $\mu_{0}%
=\cos(\theta_{0})$, $\nu_{0}=\sec^{2}(\theta_{0})$, and $\eta_{0}=0.$ If
$\theta_{0}^{p}=\pm\pi/2\mod2\pi$ take $x_{3}=\cot(\theta_{0})$, $\mu_{0}%
=\sin(\theta_{0})$, $\nu_{0}=-\csc^{2}(\theta_{0})$, and $\eta_{0}=0.$ Denote
$s_{i}=\sin(\theta_{i}-\theta_{i-1})$ and $c_{i}=\cos(\theta_{i}-\theta
_{i-1})$, for $0\leq i\leq n$. Moreover, denote by $\operatorname*{L}%
\nolimits_f\!\alpha$ the Lie derivative of a function $\alpha$ along a vector
field $f$.

Now, consider the sequence of smooth functions defined locally, for $1\leq
i\leq n$, by either
\begin{equation}%
\begin{array}
[c]{lll}%
x_{i+3} & = & \dfrac{s_{i}\nu_{i-1}+c_{i}\eta_{i-1}}{c_{i}\mu_{i-1}}\\ &  & \\
\mu_{i} & = & c_{i}\mu_{i-1}\qquad\qquad\\ \nu_{i} & = &
\operatorname*{L}\nolimits_{\tau_{1}^{i}}x_{i+3}\\
\eta_{i} & = & \operatorname*{L}\nolimits_{\tau_{2}^{i}}x_{i+3}%
\end{array}
\quad\label{ieee-regular-trailer-transformation}%
\end{equation}
when $\theta_{i}^{p}-\theta_{i-1}^{p}\neq\pm\pi/2\mod2\pi$ (\emph{regular
case}) or by
\begin{equation}%
\begin{array}
[c]{lcl}%
x_{i+3} & = & \dfrac{c_{i}\mu_{i-1}}{s_{i}\nu_{i-1}+c_{i}\eta_{i-1}}\\ &  & \\
\mu_{i} & = & s_{i}\nu_{i-1}+c_{i}\eta_{i-1}\\ \nu_{i} & = &
\operatorname*{L}\nolimits_{\tau_{1}^{i}}x_{i+3}\\
\eta_{i} & = & \operatorname*{L}\nolimits_{\tau_{2}^{i}}x_{i+3}%
\end{array}
\quad\label{ieee-singular-trailer-transformation}%
\end{equation}
when $\theta_{i}^{p}-\theta_{i-1}^{p}=\pm\pi/2\mod2\pi$ (\emph{singular case}).
It is easy to prove that, for $0\leq i\leq n$, the transformations $\phi^{i}$
defined by
$\phi^{i}(\xi_{1},\xi_{2},\theta_{0},\ldots,\theta_{i})=(x_{1},\ldots,x_{i+3})$
are smooth changes of coordinates around $p_i$ and that, moreover, we have both
$\nu_{i}(\zeta^{p_i})\neq0$ and $\mu_{i}(\zeta^{p_i})\neq0$, where~$p_i$
denotes the projection of $p$ on $\mathbb{R}^{2}\times(S^{1})^{i+1}$.

\begin{theorem}
\label{thm-trailer-conversion}For $n\geq0$, the diffeomorphism $\phi^{n}$ and
the feedback transformation $(\nu_{n},\eta_{n},\mu_{n})$ satisfy
\emph{(\ref{ieee-trailer-as-kumpera})}, and thus convert the $n$-trailer system
into a Kumpera-Ruiz normal form.
\end{theorem}

\proof We will prove that the relation (\ref{ieee-trailer-as-kumpera}) holds
for $n\geq0$ by induction on the number~$n$ of trailers. Relation
(\ref{ieee-trailer-as-kumpera}) is clearly true for $n=0$. Assume that it holds
for $n-1$ trailers. In this case we have
\[
\begin{array}[c]{lll}
\phi_{*}^{n-1}(\tau_{1}^{n-1})&=&(\nu_{n-1}\circ\psi^{n-1})\,\kappa_{1}^{n+2}\\
\phi_{*}^{n-1}(\tau_{2}^{n-1})&=&(\eta_{n-1}\circ\psi^{n-1})\,\kappa_{1}^{n+2}+
(\mu_{n-1}\circ\psi^{n-1})\,\kappa_{2}^{n+2}.
\end{array}
\]
The inductive definition of the $n$-trailer gives
\[
\begin{array}[c]{l}
\tau_{1}^{n}=\tfrac\partial{\partial\theta_{n}}\\
\tau_{2}^{n}=\sin(\theta_{n}-\theta_{n-1})\tau_{1}^{n-1}+\cos(\theta
_{n}-\theta_{n-1})\tau_{2}^{n-1}.
\end{array}
\]

Let $\phi^{n}=(\phi^{n-1},\phi_{n+3})^{T}$ be a diffeomorphism of
$\mathbb{R}^{n+3}$ such that $\phi^{n-1}$ depends on the first $n+2$
coordinates only. Let $f$ be a vector field on $\mathbb{R}^{n+3}$ of the form
$f=\alpha f^{(n-1)}+f_{n+3}$, where~$\alpha$ is a smooth function on
$\mathbb{R}^{n+3}$, the vector field $f^{(n-1)}$ is the lift of a vector field
on~$\mathbb{R}^{n+2}$ (see Notation~\ref{not-ieee-lift}), and the only non-zero
component of $f_{n+3}$ is the last one. Then we have
\begin{eqnarray}
\phi_{*}^{n}(f)&=&(\alpha\circ\psi^{n})\phi_{*}^{n-1}(f^{(n-1)})
\label{ieee-composed-derivative}+ \left(
(\mathrm{L}_{f}\phi_{n+3})\circ\psi^{n} \right) \tfrac\partial{\partial
x_{n+3}}. \nonumber
\end{eqnarray}
Note that the vector field $\phi_{*}^{n-1}(f^{(n-1)})$ is lifted along the
$x_{n+3}$-coordinate, which is defined by $\phi_{n+3}$.

In the regular case we take a regular prolongation. In this case, relation
(\ref{ieee-composed-derivative}) gives:
%$\theta_{n}(p)-\theta_{n-1}(p)\neq\pm\pi/2$ and thus:
%\[
\begin{eqnarray*}
%[c]{ccl}
\phi_{*}^{n}(\tau_{2}^{n})&=&(s_{n}\circ\psi^{n})\phi_{*}^{n-1}(\tau_{1}^{n-1})
+(c_{n}\circ\psi^{n})\phi_{*}^{n-1}(\tau_{2}^{n-1}) +\left(
(\mathrm{L}_{\tau_{2}^{n}}x_{n+3})\circ\psi^{n}\right)\tfrac\partial{\partial
x_{n+3}}\\ &=&
\left((s_{n}\nu_{n-1}+c_{n}\eta_{n-1})\circ\psi^{n}\right)\kappa_{1}^{n+2}
+\left((c_{n}\mu_{n-1})\circ\psi^{n}\right)\kappa_{2}^{n+2}+(\eta_{n}\circ\psi^{n})\kappa_{1}^{n+3}\\
&=& (\eta_{n}\circ\psi^{n})\kappa_{1}^{n+3}+\left(
c_{n}\mu_{n-1}\circ\psi^{n}\right) \times \left[  \left(
\tfrac{s_{n}\nu_{n-1}+c_{n}\eta_{n-1}}{c_{n}\mu_{n-1}}
\circ\psi^{n}\right)\kappa_{1}^{n+2}+\kappa_{2}^{n+2}\right] \\
&=&(\eta_{n}\circ\psi^{n})\kappa_{1}^{n+3} +(\mu_{n}\circ\psi^{n})\left(
x_{n+3}\kappa_{1}^{n+2}+\kappa_{2}^{n+2}\right)\\
&=&(\eta_{n}\circ\psi^{n})\kappa_{1}^{n+3}+(\mu_{n}\circ\psi^{n})\kappa_{2}^{n+3}.
\end{eqnarray*}
%\]
In the singular case we take the singular prolongation and,
by~(\ref{ieee-composed-derivative}), we obtain:
%\[
\begin{eqnarray*}
%[c]{ccl}
\phi_{*}^{n}(\tau_{2}^{n}) & = & \left(
(s_{n}\nu_{n-1}+c_{n}\eta_{n-1})\circ\psi^{n}\right) \times \left[
\kappa_{1}^{n+2}+\left(  \tfrac{c_{n}\mu_{n-1}}{s_{n}\nu_{n-1}+
     c_{n}\eta_{n-1}}\circ\psi^{n}\right)\kappa_{2}^{n+2}\right]\\
& & \mbox{} +  (\eta_{n}\circ\psi^{n})\kappa_{1}^{n+3}\\
&=&(\eta_{n}\circ\psi^{n})\kappa_{1}^{n+3} +(\mu_{n}\circ\psi^{n})\left(
\kappa_{1}^{n+2}+x_{n+3}\kappa_{2}^{n+2}\right)\\
&=&(\eta_{n}\circ\psi^{n})\kappa_{1}^{n+3}+(\mu_{n}\circ\psi^{n})\kappa_{2}^{n+3}.
\end{eqnarray*}
%\]
Moreover, in both cases, we have
%\[
%\begin{array}
%[c]{ccl}
%\phi_{*}^{n}(\tau_{1}^{n}) & = & \left(  (\mathrm{L}_{\tau_{1}^{n}}%
%x_{n+3})\circ\psi^{n}\right)  \tfrac\partial{\partial x_{n+3}}\\
%& = & (\nu_{n}\circ\psi^{n})\kappa_{1}^{n+3}.%
%\end{array}
%\]
\[
\phi_{*}^{n}(\tau_{1}^{n}) = \left(
(\mathrm{L}_{\tau_{1}^{n}}x_{n+3})\circ\psi^{n}\right) \tfrac\partial{\partial
x_{n+3}}=(\nu_{n}\circ\psi^{n})\kappa_{1}^{n+3}.
\]
It follows that, both in the regular and in the singular case, relation
(\ref{ieee-trailer-as-kumpera}) holds for $n\geq0$. \endproof

\medskip

Reversing the construction given in the previous proof leads to the following
surprising result \cite{pasillas-respondek-nolcos}, which states that the
$n$-trailer system is a universal local model for all Goursat structures
--- see also \cite{montgomery-zhitomirskii} and
\cite{pasillas-respondek-goursat}.

\begin{theorem}
\label{ieee-thm-trailer-universal}Any Goursat structure on a manifold $M$ of
dimension $n+3$ is equivalent, in a small enough neighborhood of any point $q$
in $M$, to the $n$-trailer considered around a suitably chosen point $p$ of its
configuration space $\mathbb{R}^{2}\times(S^{1})^{n+1}$.
\end{theorem}

\proof By Theorem~\ref{ieee-thm-kumpera-ruiz}, our Goursat structure is, in a
small enough neighborhood of any point $q\ $in $M$, equivalent to a
Kumpera-Ruiz normal form $\kappa^{n+3}$. Denote by $y=(y_{1},\ldots,y_{n+3})$
the coordinates of $\kappa^{n+3}$ and put
$(y_{1}^{q},\ldots,y_{n+3}^{q})=y(q)$.

Recall that, by definition, the pair of vector fields $\kappa^{n+3}$ is given
by a sequence of prolongations $\kappa^{i}=\sigma_{i-3}\circ\cdots\circ\sigma
_{1}(\kappa^{3})$, where $\sigma_{j}$ belongs to $\{R_{c_{j}},S\}$, for $1\leq
j\leq i-3$ and $3\leq i\leq n+3$. We call a coordinate $y_{i}$ such that
$\kappa^{i}=S(\kappa^{i-1})$ a \emph{singular coordinate}, and a
coordinate~$y_{i}$ such that $\kappa^{i}=R_{c}(\kappa^{i-1})$ a \emph{regular
coordinate}. It follows from the proof of Theorem~\ref{ieee-thm-kumpera-ruiz}
--- see \cite{cheaito-mormul} and \cite{pasillas-respondek-goursat} --- that for
all singular coordinates we have $y_{i}^{q}=0$; but for regular coordinates,
the constants $y_{i}^{q}$ can be arbitrary real numbers.

To prove the theorem, we will define a point $p$ of $\mathbb{R}%
^{2}\times(S^{1})^{n+1}$ whose coordinates $\zeta(p)=(\xi_{1}^{p},\xi_{1}%
^{p},\theta_{0}^{p},\ldots,\theta_{n}^{p})$ satisfy $(x\circ\zeta)(p)=y(q)$,
where $x$ and $\zeta$ denote the coordinates used in the Proof of Theorem
\ref{thm-trailer-conversion}. First, put the axle of the last trailer at
$(y_{2}^{q},y_{1}^{q})$ and take $\theta_{0}^{p}=\arctan(y_{3}^{q})$. Compute
 $x_{3}=\tan(\theta_0)$, $\mu_{0}=\cos(\theta_0)$, $\nu_{0}=\sec^2(\theta_0)$, and
$\eta_{0}=0$. Then, take for $i=1$ up to $n$, the following values for the
angles $\theta_{i}^{p}\mod 2\pi$. If the
coordinate $y_{i+3}$ is singular then put $\theta_{i}^{p}=\theta_{i-1}^{p}%
+\pi/2$ and compute the coordinate $x_{i+3}$ and the smooth functions $\mu
_{i}$, $\nu_{i}$, and $\eta_{i}$ using~(\ref{ieee-singular-trailer-transformation}%
). If $y_{i+3}$ is regular then put
\[
\theta_{i}^{p}=\arctan\left(  \dfrac{\mu_{i-1}(p)y_{i+3}^{q}-\eta_{i-1}%
(p)}{\nu_{i-1}(p)}\right)  +\theta_{i-1}^{p}
\]
and compute the coordinate $x_{i+3}$ and the smooth functions $\mu_{i}$,
$\nu_{i}$, and $\eta_{i}$ using~(\ref{ieee-regular-trailer-transformation}).
%The result of this construction is that $(x\circ\zeta)(p)=y(q)$.
By Theorem \ref{thm-trailer-conversion}, the coordinates $x\circ\zeta$ convert
the $n$-trailer into a Kumpera-Ruiz normal form. By the above defined
construction, this normal form has the same singularities as~$\kappa^{n+3}$ and
is defined around the same point of $\mathbb{R}^{n+3}$ (if we translate the
coordinates in order to center them then these Kumpera-Ruiz normal forms have
the same constants in the regular prolongations). Hence, the diffeomorphism
$\zeta^{-1} \circ x^{-1} \circ y$ gives the claimed equivalence.
\endproof

%Since, by Theorem \ref{thm-trailer-conversion}, the coordinates $x\circ\zeta$
%and $y$ convert the $n$-trailer and our Goursat structure, respectively, into
%Kumpera-Ruiz normal forms that have the same singularities and that are defined
%around the same points, our Goursat structure at~$q$ is locally equivalent to
%the~$n$-trailer system at~$p$.\hfill$\square$

\section{Motion Planning}

\label{sec-motion-planning}

As we have proved in Section~\ref{sec-nilpotentization}, the Lie algebras
generated by Kumpera-Ruiz normal forms are nilpotent. This property is
fundamental \citeaffixed{laferriere-sussmann,liu-approximation}{see e.g.}
because it allows to solve the nonholonomic motion planning problem in the case
of Goursat structures. But Kumpera-Ruiz normal forms have also an other
interesting property: they are ``triangular''  --- in the sense of
\citeasnoun{murray-sastry}, see also \cite{marigo}. Indeed, it follows directly
from their construction that they give a control system $\dot x=\kappa_{1}
^{n}(x)\,u_{1}+\kappa_{2}^{n}(x)\,u_{2}$ that can be written --- see
\cite{cheaito-mormul}, \cite{cheaito-mormul-pasillas-respondek},
\cite{kumpera-ruiz}, and \cite{pasillas-respondek-goursat}, after a permutation
of the $x_{i}$'s, in the following form:
\[%
\begin{array}
[c]{lcl}%
\dot x_{1} & = & u_{1}\\ \dot x_{2} & = & u_{2}\\ \dot x_{3} & = &
f_{3}(x_{1},x_{2})\,u_{2}\\ & \vdots & \\ \dot x_{n-1} & = &
f_{n-3}(x_{1},\ldots,x_{n-2})\,u_{2}\\ \dot x_{n} & = &
f_{n}(x_{1},\ldots,x_{n-2},x_{n-1})\,u_{2},
\end{array}
\]
where, for $3\leq i\leq n$, the $f_{i}$'s are polynomials. Therefore, if we
take polynomial controls of the form
\begin{equation}%
\begin{array}
[c]{l}%
u_{1}(t)=a_{0}+a_{1}t+\cdots+a_{n-2}t^{n-2}\\ u_{2}(t)=b_{0},
\end{array}
\label{nice-controls}%
\end{equation}
as it has been proposed in \cite{tilbury-murray-sastry} in order to steer the
$n$-trailer around regular points, then the control system $\dot x=\kappa_{1}%
^{n}(x)\,u_{1}+\kappa_{2}^{n}(x)\,u_{2}$ can be integrated by successive
quadratures, which leads to a system of $n$ polynomial equations
\begin{equation}
x_{i}^{T} = P_{i}(a_{0},\ldots,a_{n-2},b_{0},x_{1}^{0},\ldots,x_{n}^{0}),
\label{polynomial-system}
\end{equation}
%\begin{equation}%
%\begin{array}
%[c]{ccc}%
%x_{1}^{T} & = & P_{1}(a_{0},\ldots,a_{n-1},b_{0},x_{1}^{0},\ldots,x_{n}^{0})\\
%& \vdots & \\
%x_{n}^{T} & = & P_{n}(a_{0},\ldots,a_{n-1},b_{0},x_{1}^{0},\ldots,x_{n}^{0}),
%\end{array}
%\label{polynomial-system}%
%\end{equation}
where, for $1\leq i\leq n$, the $P_{i}$'s are polynomial functions of all their
arguments and the vectors $x^{0}$ and $x^{T}$ denote the initial and final
configuration, respectively. Hence, for Goursat structures, the nonholonomic
motion planning problem can be transformed into an algebraic problem: the
resolution of a system of polynomial equations.

In the regular case, the polynomial system (\ref{polynomial-system}) is
actually a linear system, which reflects the flatness of the system around
regular configurations
\citeaffixed{fliess-intro-flat,martin-rouchon-driftless}{see e.g.}. Indeed, if
we fix $b_{0}\neq0$ then the final condition $x^{T}$ belongs to a hyperplane
$E_{b_{0}}$ of~$\mathbb{R}^{n}$, parameterized by the constants
$a_{0},...,a_{n-2}$, and the motion planning problem leads to the resolution of
a full-rank system of linear equations. Note, however, that this approach fails
if $b_{0}=0$. That is, when the control~(\ref{nice-controls}) produces an
abnormal trajectory.

In the singular case, the system (\ref{polynomial-system}) is truly polynomial,
which causes at least two problems. Firstly, in general, it cannot be solved
explicitly and we must use numerical approximations. Secondly, although the
system is globally controllable, there may be restrictions on the set of points
accessible with controls of the form (\ref{nice-controls}) and, in general,
describing these restrictions leads to another system of polynomial equations.

Fortunately, in small dimension the situation is not so complicated. Indeed, if
we consider a mobile robot towing two or three trailers then
(\ref{polynomial-system}) leads to a system of linear equations together with a
single quadratic equation, which can be solved easily, and therefore the
nonholonomic motion planning problem admits an explicit solution. Moreover, the
set of points accessible by the family of controls (\ref{nice-controls}) can be
completely characterized: if we fix all initial and final angles, the set of
points in the $(\xi_{1},\xi_{2})$-plane that can be reached from a given
configuration is delimited by a parabola.

\section{The Two-Trailer System}

\label{sec-two-trailer}

The results presented in this paper can be directly applied to the mobile robot
towing two trailers, which is given by
\[%
\begin{array}
[c]{l}%
\dot\xi_{1}=\cos(\theta_{2}-\theta_{1})\cos(\theta_{1}-\theta_{0})\cos
(\theta_{0})v_{2}\\
\dot\xi_{2}=\cos(\theta_{2}-\theta_{1})\cos(\theta_{1}-\theta_{0})\sin
(\theta_{0})v_{2}\\
\dot\theta_{0}=\cos(\theta_{2}-\theta_{1})\sin(\theta_{1}-\theta_{0})v_{2}\\
\dot\theta_{1}=\sin(\theta_{2}-\theta_{1})v_{2}\\ \dot\theta_{2}=v_{1}.
\end{array}
\]
We will consider this system in a small enough neighborhood of the singular
locus $\theta_{2}-\theta_{1}=\pm\frac\pi2$, where the transformation
\begin{equation}
\label{two-trailer-transformation}
\begin{array}[c]{l}
x_{1} =\xi_{1}\\ x_{2} =\xi_{2}\\ x_{3} =\tan(\theta_{0})\\ x_{4}
=\sec(\theta_{0})^{3}\tan(\theta_{1}-\theta_{0})\\ x_{5}
=\tfrac{\cos(\theta_{0})^{4}\cos(\theta_{1}-\theta_{0})\cos
(\theta_{2}-\theta_{1})}{3\tan(\theta_{0})\tan(\theta_{1}-\theta_{0}%
)\sin(\theta_{1}-\theta_{0})\cos(\theta_{2}-\theta_{1})+\sec(\theta_{1}%
-\theta_{0})^{2}\left(  \sin(\theta_{2}-\theta_{1})-\sin(\theta_{1}-\theta
_{0})\cos(\theta_{2}-\theta_{1})\right)  }%
\end{array}
\end{equation}
converts the two-trailer into the following Kumpera-Ruiz normal form:
\[%
\begin{array}
[c]{l}%
\dot x_{1}=x_{5}\,u_{2}\\ \dot x_{2}=x_{5}x_{3}\,u_{2}\\ \dot
x_{3}=x_{5}x_{4}\,u_{2}\\ \dot x_{4}=u_{2}\\ \dot x_{5}=u_{1}.
\end{array}
\]
The feedback transformation between the old and new controls is given in the
Matlab file \verb|system.m| of \citeasnoun{pasillas-thesis}, which is available
by e-mail request to the first author.

\bigskip

We will show that (\ref{two-trailer-transformation}) defines a diffeomorphism
of a well chosen set $V$ onto ${\mathbb R}^{5}$.  In order to do so, put
$\tilde \theta_{0}=\theta_{0}$, $\tilde \theta_{1}=\theta_{1}-\theta_{0}$,
$\tilde \theta_{2}=\theta_{2}-\theta_{1}$ and represent $x_{5}$, defined by
(\ref{two-trailer-transformation}), as
\begin{equation}
x_{5}=\frac{ a\cos\tilde\theta_{2} }{ b\cos\tilde\theta_{2}+c\sin
\tilde\theta_{2} }, \label{ijc-thesis-map}
\end{equation}
where $a$, $b$, and $c$ are smooth functions of $\tilde\theta_{0}$ and
$\tilde\theta_{1}$. It is easy to observe that for any fixed values of
$\tilde\theta_{0}$ and $\tilde\theta_{1}$ and on any interval
$(-\pi/2+k\pi,\pi/2+k\pi)$ the function $b\cos\tilde\theta_{2}+c\sin
\tilde\theta_{2}$ vanishes exactly one time. Let $\gamma(\tilde\theta_{0},
\tilde\theta_{1})
\in
(-\pi/2+k\pi,\pi/2+k\pi)$ and $\delta(\tilde\theta_{0}, \tilde\theta_{1}) \in
(-\pi/2+(k+1)\pi,\pi/2+(k+1)\pi)$ be zeros of $b\cos\tilde\theta_{2}+c\sin
\tilde\theta_{2}$. We claim  that for any fixed values of $\tilde\theta_{0}$
and $\tilde\theta_{1}$, the map(\ref{ijc-thesis-map}) establishes a
diffeomorphism of the interval $(\gamma(\tilde\theta_{0}, \tilde\theta_{1}),
\delta(\tilde\theta_{0}, \tilde\theta_{1}))$ onto ${\mathbb R}$. To see this,
firstly, observe that $\delta(\tilde\theta_{0},
\tilde\theta_{1})-\gamma(\tilde\theta_{0}, \tilde\theta_{1})=\pi$. Secondly,
dividing, on well chosen sets,  the numerator and denominator of the right hand
side of (\ref{ijc-thesis-map}) by $\sin\tilde\theta_{2}$ and replacing
$\cot\tilde\theta_{2}$ by $-y$ or  by $\cos\tilde\theta_{2}$ and replacing
$\tan\tilde\theta_{2}$ by $y$, we represent $x_{5}$, respectively, as
restrictions of either the homography $x_{5}=\frac{-ay}{-by+c}$ or
$x_{5}=\frac{a}{cy+b}$. In both cases, derivative of $x_{5}$ with respect to
$\tilde\theta_{2}$ is of the same sign which is the sign of $-ac$. This and the
fact that $b\cos\tilde\theta_{2}+c\sin \tilde\theta_{2}$ vanishes at
$\gamma(\tilde\theta_{0}, \tilde\theta_{1})$ and $\delta(\tilde\theta_{0},
\tilde\theta_{1})$ implies that when $\tilde\theta_{2}$ changes between
$\gamma(\tilde\theta_{0}, \tilde\theta_{1})$ and $\delta(\tilde\theta_{0},
\tilde\theta_{1})$ then the value of $x_{5}$ either grows monotonically between
minus and plus infinity or decreases monotonically between  plus and minus
infinity.

Now define $V=\{(\xi_{1},\xi_{2}, \tilde \theta_{0}, \tilde \theta_{1}, \tilde
\theta_{2}): \xi_{1}\in \mathbb{R},\ \xi_{2}\in \mathbb{R},\ \tilde
\theta_{0}\in(-\pi/2,\pi/2),\ \tilde \theta_{1}\in(-\pi/2,\pi/2),\ \tilde
\theta_{2}\in(\gamma(\tilde\theta_{0},
\tilde\theta_{1}),\delta(\tilde\theta_{0}, \tilde\theta_{1}))\}.$ The above
considerations imply that (\ref{two-trailer-transformation}) establishes a
diffeomorphism between~$V$ and~${\mathbb R}^{5}$.

\bigskip

Now fix an initial condition $\xi^{0}=(\xi_{1}^{0},\xi_{2}^{0},\theta_{0}^{0},
\theta_{1}^{0},  \theta_{2}^{0})$ and a terminal condition
$\xi^{T}=(\xi_{1}^{T},\xi_{2}^{T},\theta_{0}^{T},\theta_{1}^{T},\theta_{2}^{T}).$
Let $k$ be  an integer  such that
$\tilde\theta_{2}^{0}\in(\gamma(\tilde\theta_{0}^{0}, \tilde\theta_{1}^{0}),
\delta(\tilde\theta_{0}^{0}, \tilde\theta_{1}^{0})),$ where
$\gamma(\tilde\theta_{0}^{0}, \tilde\theta_{1}^{0})
\in
(-\pi/2+k\pi,\pi/2+k\pi)$. Assume that $\gamma(\tilde\theta_{0}^{T},
\tilde\theta_{1}^{T})$ satisfies $\gamma(\tilde\theta_{0}^{T},
\tilde\theta_{1}^{T})\in (-\pi/2+k\pi,\pi/2+k\pi)$ for the same $k$, where
$\tilde\theta_{2}^{T}\in(\gamma(\tilde\theta_{0}^{T}, \tilde\theta_{1}^{T}),
\delta(\tilde\theta_{0}^{T}, \tilde\theta_{1}^{T})$. Then any smooth curve
$x(t)$ in~${\mathbb R}^{5}$ such that $x(0)=\xi^{0}$ and $x(T)=\xi^{T}$ has a
unique smooth preimage, with respect to (\ref{two-trailer-transformation}).

\bigskip

The controls $u_{1}$ and $u_{2}$ that steer the above Kumpera-Ruiz normal form
from $p=(p_{1},p_{2},p_{3},p_{4},p_{5})$ at $t=0$ to $q=(q_{1},q_{2}%
,q_{3},q_{4},q_{5})$ at $t=1$ are given in \citeasnoun{pasillas-thesis} (see
the files \verb|main.m| and \verb|system.m|). Here, we will consider the
particular case $p=(0,0,0,0,0)$, which is simpler to analyze. If we take a pair
of polynomial controls of the form
\begin{equation}%
\begin{array}
[c]{l}%
u_{1}(t)=a_{2}+a_{3}t+a_{4}t^2+a_{5}t^3\\ u_{2}(t)=a_{1},
\end{array}
\label{matlab-controls}%
\end{equation}
then we will obtain:
\[%
\begin{array}
[c]{ccl}%
x_{1}(1) & = & a_{1}\left(  \frac{1}{2}a_{2}+\frac{1}{6}a_{3}+\frac{1}%
{12}a_{4}+\frac{1}{20}a_{5}\right)  \\
x_{2}(1) & = & a_{1}^{3}\left(  \frac{1}{15}a_{2}^{2}+\frac{1}{112}a_{3}%
^{2}+\frac{1}{405}a_{4}^{2}+\frac{1}{1056}a_{5}^{2}\right.  \\
&  & \left.  +\frac{7}{144}a_{2}a_{3}+\frac{8}{315}a_{2}a_{4}+\frac{5}%
{320}a_{2}a_{5}\right.  \\
&  & \left.  +\frac{3}{320}a_{3}a_{4}+\frac{5}{864}a_{3}a_{5}+\frac{11}%
{3600}a_{4}a_{5}\right)  \\
x_{3}(1) & = & a_{1}^{2}\left(  \frac{1}{3}a_{2}+\frac{1}{8}a_{3}+\frac{1}%
{15}a_{4}+\frac{1}{24}a_{5}\right)  \\ x_{4}(1) & = & a_{1}\\ x_{5}(1) & = &
a_{2}+\frac{1}{2}a_{3}+\frac{1}{3}a_{4}+\frac{1}{4}a_{5}.
\end{array}
\]
Now, if we want the final point to be $q$, we should take:
\[%
\begin{array}
[c]{l}%
a_{1}=q_{4}\\
a_{3}=-12a_{2}-\frac{360}{q_{4}^{2}}q_{3}+\frac{240}{q_{4}}q_{1}+12q_{5}\\
a_{4}=30a_{2}+\frac{1440}{q_{4}^{2}}q_{3}-\frac{900}{q_{4}}q_{1}-60q_{5}\\
a_{5}=-20a_{2}-\frac{1200}{q_{4}^{2}}q_{3}+\frac{720}{q_{4}}q_{1}+60q_{5}.
\end{array}
\]
The constraint $x_{2}(1)=q_{2}$ gives, for $a_{2}$, the quadratic equation
\[%
\begin{array}
[c]{r}%
\frac{1}{55440}q_{4}^{3}\;a_{2}^{2}+\left(  \frac{1}{462}q_{3}q_{4}-\frac
{2}{693}q_{1}q_{4}^{2}+\frac{1}{11088}q_{5}q_{4}^{3}\right)  \;a_{2}%
-\frac{179}{462}q_{1}q_{3}\text{\qquad}\\
+\frac{60}{77}\frac{q_{3}^{2}}{q_{4}}+\frac{15}{77}q_{1}^{2}q_{4}-\frac
{1}{154}q_{3}q_{4}q_{5}-\frac{1}{308}q_{1}q_{4}^{2}q_{5}+\frac{1}{2772}%
q_{4}^{3}q_{5}^{2}=q_{2}.
\end{array}
\]
It follows that a point $q$ is reachable from zero, with controls of the form
(\ref{matlab-controls}), if and only if:
\begin{equation}
\label{parabola}
\begin{array}[c]{r}
q_{2}\geq\frac{5}{63}q_{4}\;q_{1}^{2}+\left(  -\frac{3}{14}q_{3}+\frac{1}%
{252}q_{4}^{2}q_{5}\right)  \;q_{1}+\frac{5}{7}\frac{q_{3}^{2}}{q_{4}}%
-\frac{1}{84}q_{3}q_{4}q_{5}+\frac{1}{4032}q_{4}^{3}q_{5}^{2}.
\end{array}
\end{equation}

To steer the two-trailer system: (i) Fix initial and terminal configurations
$\zeta^{0}$ and $\zeta^{T}$ that can be joined by a control of the
form~(\ref{matlab-controls}). In the particular case $x(\zeta^{0})=0$, the
terminal configuration must satisfy (\ref{parabola}). Check that $\zeta^{0}$
and $\zeta^{T}$ belong to the domain~$V$ defined above. (ii) Compute, using the
transformation~(\ref{two-trailer-transformation}), the initial and final
configurations $x^{0}=x(\zeta^{0})$ and $x^{T}=x(\zeta^{T})$ in the coordinates
of the Kumpera-Ruiz normal form. (iii) Compute a pair of controls $(u_1,u_2)$
of the form~(\ref{matlab-controls}), that steer the system from $x^{0}$ to
$x^{T}$. (iv) Apply to the original system the controls
$v_{1}=\widehat{\nu}_{5}\,u_{1}+\widehat{\eta}_{5}\,u_{2}$ and
$v_{2}=\widehat{\mu}_{5}\,u_{2}$, given
by~(\ref{inverse-feedback-transformation}). For additional details, we refer
the reader to the three Matlab files \verb|main.m|, \verb|system.m|, and
\verb|chgcoord.m| of \citeasnoun{pasillas-thesis}.

We plot, in Figure~1, the trajectory joining the points
$\zeta^{0}=(0,0,0,0,\tfrac\pi4)$ and
$\zeta^{T}=(0,1,0,\tfrac\pi4,\tfrac{3\pi}4)$; and, in Figure~2, the trajectory
joining the points $\zeta^{0}=(0,0,0,-\tfrac\pi4,0)$ and
$\zeta^{T}=(0,1,0,\tfrac\pi4,\tfrac{3\pi}4)$. In both cases, illustrated by
Figures 1 and 2, the initial and terminal configurations satisfy the above
described conditions which justifies the application of our method. Indeed, in
the first case we have for the initial configuration $\tilde
\beta_{2}^{0}=\pi/4$ so
 $\tilde \beta^{0}\in (\gamma(0,0),
\delta(0,0))$, where $\gamma(0,0)=0$, $\delta(0,0)=\pi$. For the terminal
configuration we have $\tilde \beta_{2}^{T}=\pi/2\in \gamma(0,\pi/4),
\delta(0,\pi/4)$, where $\gamma(0,\pi/4)=\arctan(\sqrt{2}/4)$ which obviously
satisfies the required condition. Similar calculations hold in the case
considered in Figure 2. Observe that, in both cases, the terminal configuration
$\zeta^{T}$ belongs to the singular locus.

\newpage

\begin{figure}[ht]
\begin{minipage}[h]{0.425\hsize}
\begin{center}
\includegraphics[width=\hsize]{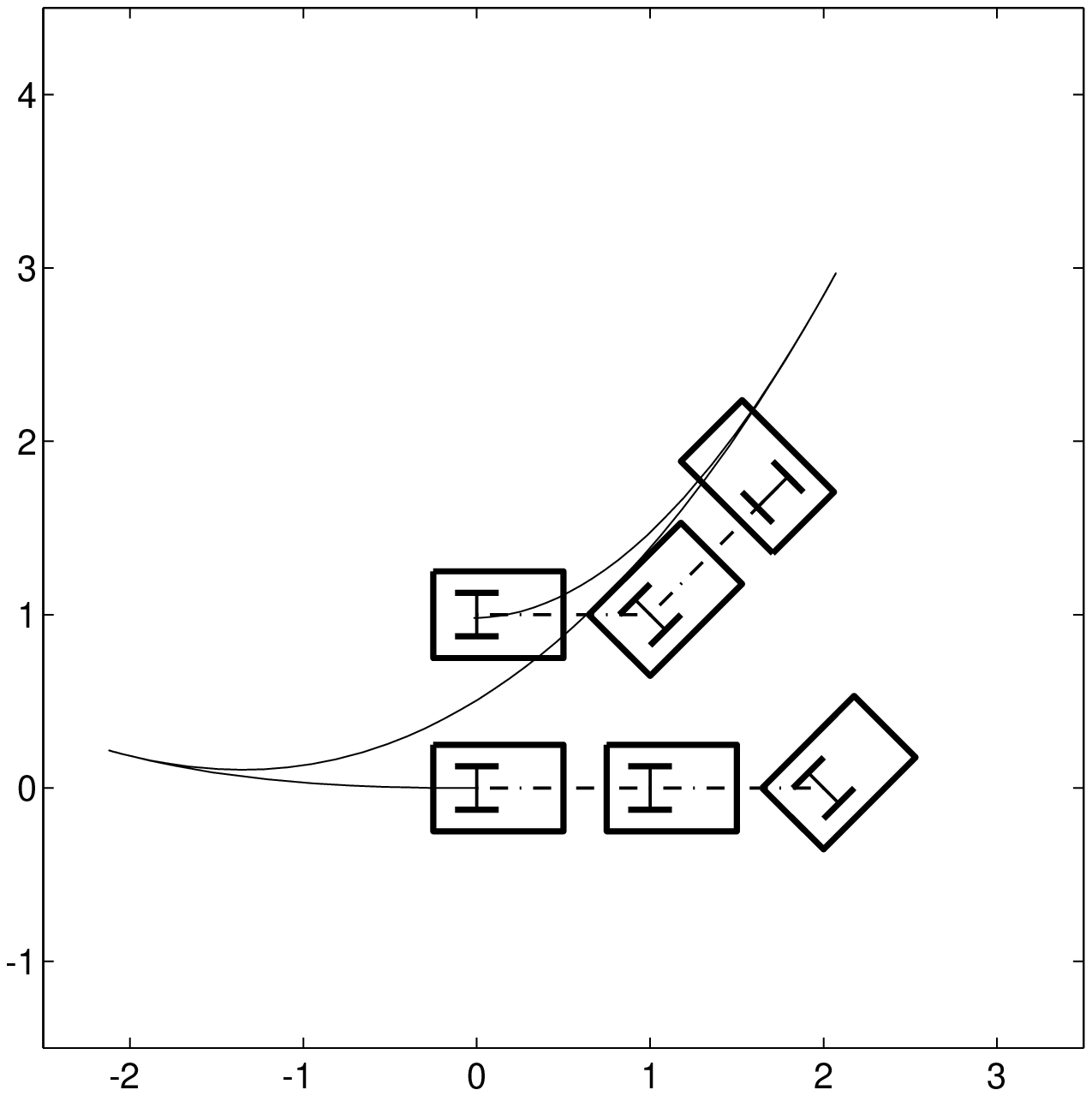}
\end{center}
\end{minipage}
\hspace{0.1\hsize}
\begin{minipage}[h]{0.425\hsize}
\begin{center}
\includegraphics[width=\hsize]{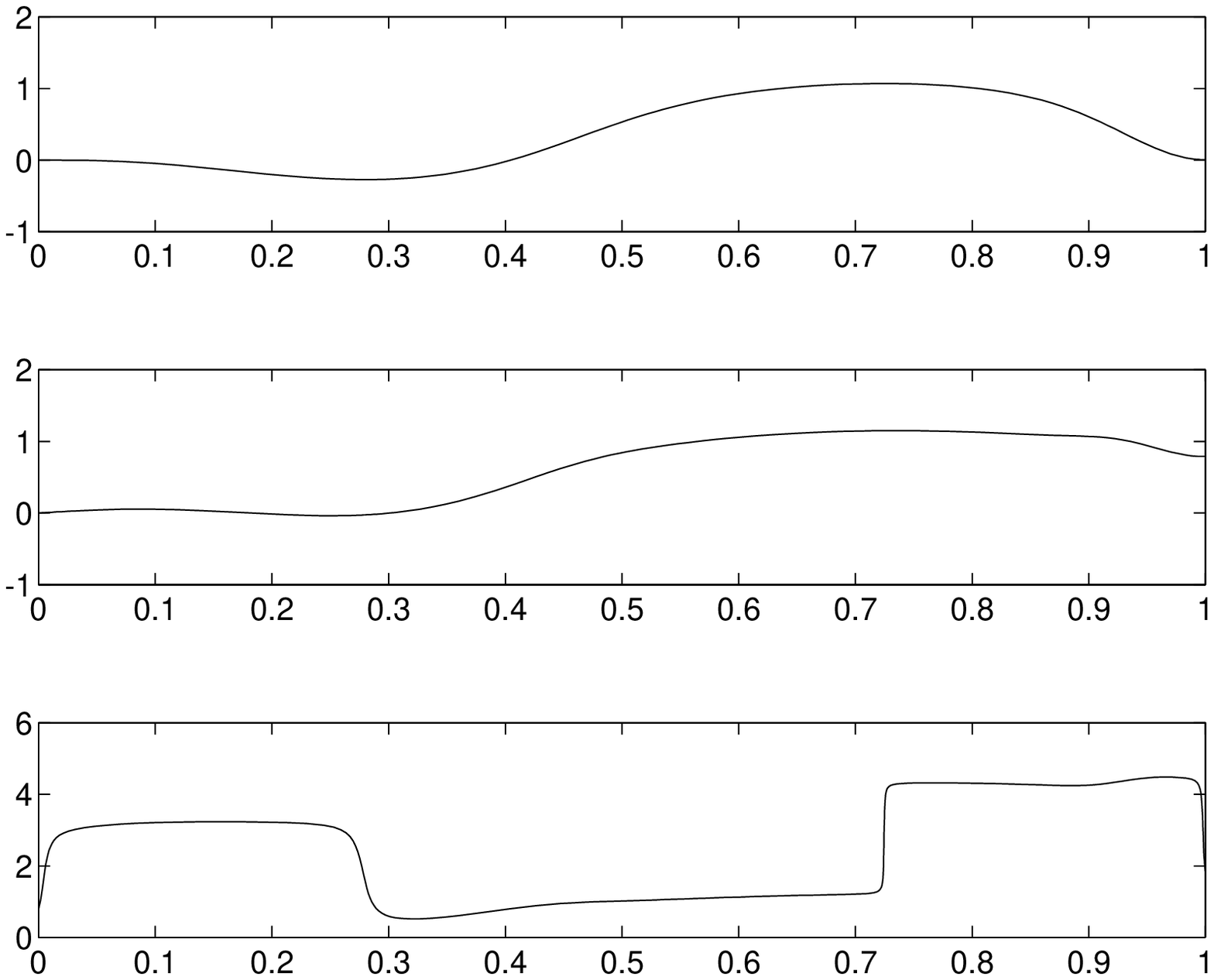}
\end{center}
\end{minipage}
\caption{(a) Initial and final configurations. Trajectory of the last trailer.
         (b) From the top: the angles $\theta_{0}$, $\theta_{1}$ and $\theta_{2}$.}
\end{figure}

\bigskip

\begin{figure}[ht]
\begin{minipage}[h]{0.425\hsize}
\begin{center}
\includegraphics[width=\hsize]{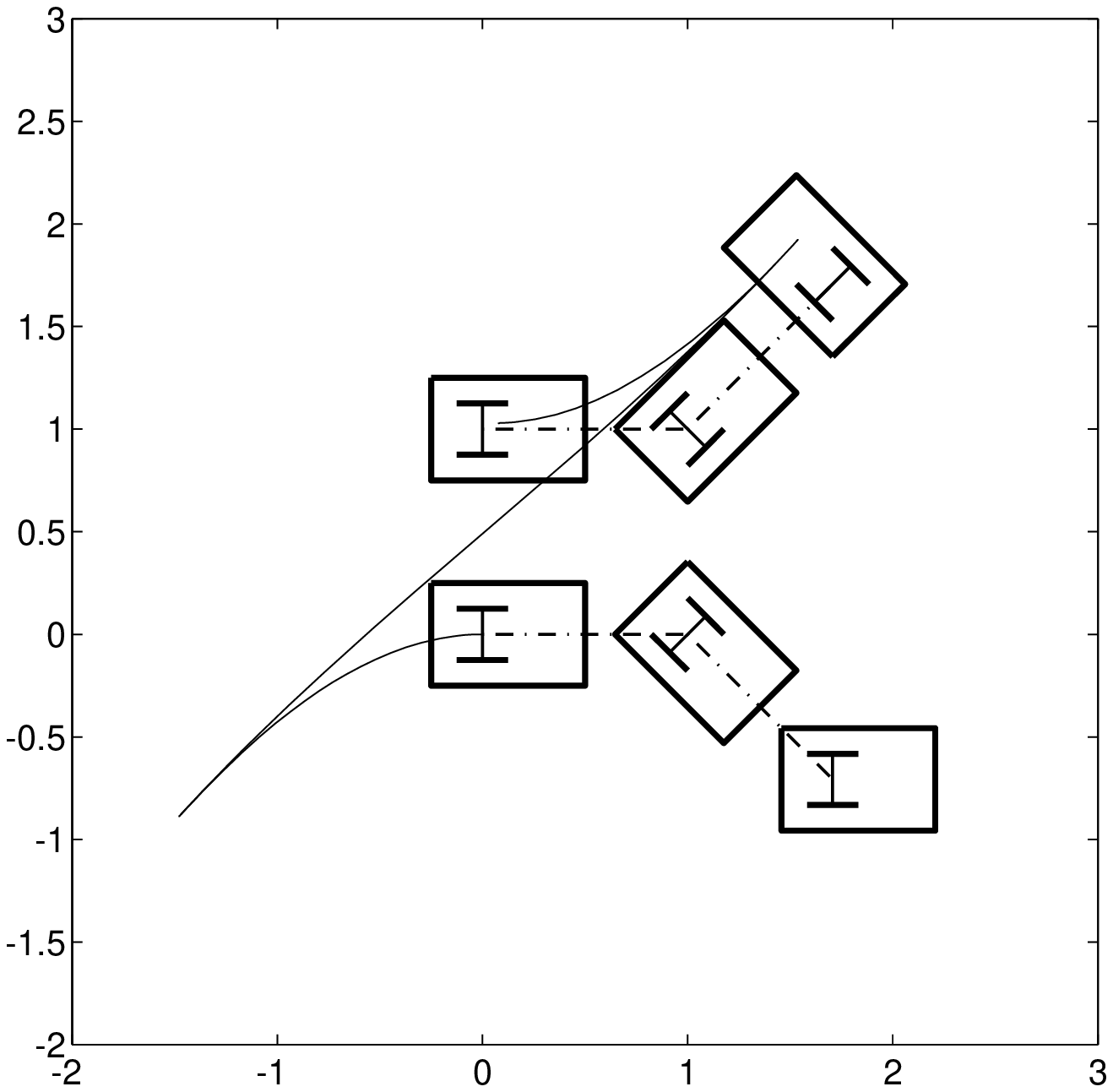}
\end{center}
\end{minipage}
\hspace{0.1\hsize}
\begin{minipage}[h]{0.425\hsize}
\begin{center}
\includegraphics[width=\hsize]{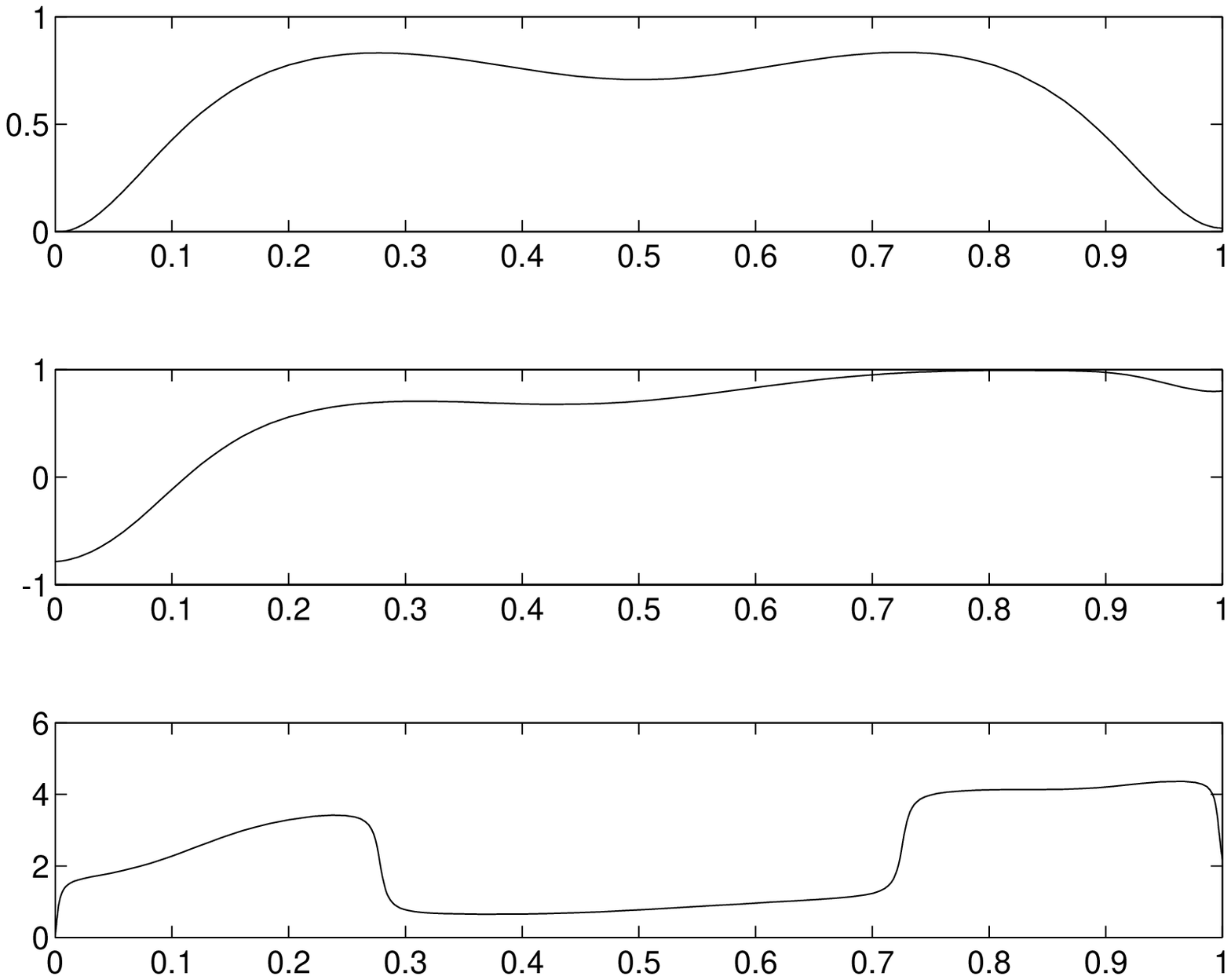}
\end{center}
\end{minipage}
\caption{(a) Initial and final configurations. Trajectory of the last trailer.
         (b) From the top: the angles $\theta_{0}$, $\theta_{1}$ and $\theta_{2}$.}
\end{figure}

\end{document}